\documentclass{svproc}
\usepackage{url}

\usepackage{amssymb,amsmath,amsfonts}
\usepackage{mathrsfs}

\def\Z{\mathbb{Z}}

\begin{document}
\mainmatter              % start of a contribution
\title{Uniform Approximation by Polynomials with Integer Coefficients via the Bernstein Lattice}
\titlerunning{The Bernstein Lattice}  % abbreviated title (for running head)
\author{C.~Sinan G\"unt\"urk\inst{1} \and Weilin Li\inst{2}}
\authorrunning{G\"unt\"urk and Li} % abbreviated author list (for running head)
\institute{NYU Courant Institute, New York, NY, USA\\
\email{gunturk@courant.nyu.edu},\\ 
\and
CUNY City College, New York, NY, USA \\
\email{wli6@ccny.cuny.edu}}

\maketitle              % typeset the title of the contribution

\begin{abstract}
Let $\mathscr{C}_\mathbb{Z}([0,1])$ be the metric space of real-valued continuous functions on $[0,1]$ with integer values at $0$ and $1$, equipped with the uniform (supremum) metric $d_\infty$. It is a classical theorem in approximation theory that the ring $\mathbb{Z}[X]$ of polynomials with integer coefficients, when  considered as a set of functions on $[0,1]$, is dense in $\mathscr{C}_\mathbb{Z}([0,1])$. In this paper, we offer a strengthening of this result by identifying a substantially small subset $\bigcup_n \mathscr{B}_n$ of $\mathbb{Z}[X]$ which is still dense in $\mathscr{C}_\mathbb{Z}([0,1])$. Here $\mathscr{B}_n$, which we call the ``Bernstein lattice,'' is the lattice generated by the polynomials 
$$p_{n,k}(x) := \binom{n}{k} x^k(1-x)^{n-k}, ~~k=0,\dots,n.$$
Quantitatively, we show that for any $f \in \mathscr{C}_\mathbb{Z}([0,1])$,
$$d_\infty(f, \mathscr{B}_n) \leq \frac{9}{4} \omega_f(n^{-1/3}) + 2 n^{-1/3}, ~~n \geq 1,$$
where $\omega_f$ stands for the modulus of continuity of $f$. We also offer a more general bound which can be optimized to yield better decay of approximation error for specific classes of continuous functions.

% We would like to encourage you to list your keywords within
% the abstract section using the \keywords{...} command.
\end{abstract}

\section{Introduction and the statement of the Main Theorem}

The topic of approximation using polynomials with integer coefficients has a mature history going back to the early 20th century. The earliest known result appears to be by P{\'a}l \cite{pal1914} who showed in 1914 that for any $0 < \alpha < 1$, any $f \in \mathscr{C}([-\alpha,\alpha])$ can be approximated uniformly by elements of $\Z[X]$, provided that $f(0) \in \Z$. Here  $\mathscr{C}([a,b])$ is the metric space of continuous real-valued functions on the interval $[a,b]$, equipped with the uniform (supremum) distance 
$$d_\infty(f,g):= \|f - g \|_{[a,b]} := \sup_{[a,b]} |f - g|,$$
and
$\Z[X]$ stands for the commutative ring of polynomials with integer coefficients which we consider as a subset of $\mathscr{C}([a,b])$.

The sufficient condition of P{\'a}l is clearly also necessary since any polynomial with integer coefficients is integer-valued at $0$, which is a property preserved in the limit. Hence, for $\alpha < 1$, the conditions $f \in \mathscr{C}([-\alpha,\alpha])$ and $f(0)\in \Z$ together characterize uniform approximability on $[-\alpha, \alpha]$ by elements of $\Z[X]$.

Kakeya \cite{kakeya1914} extended P{\'a}l's result to the space $\mathscr{C}([-1,1])$ via the necessary and sufficient condition $$f(-1),f(0),f(1),\frac{f(-1)+f(1)}{2} \in \Z.$$
The fourth of these arithmetic conditions is easily checked to be satisfied by elements of $\Z[X]$ and their limits.

Other early important results of the theory developed in the 1920s and 30s include those by Fekete, Chlodovsky, Bernstein and Kantorovich (see \cite{chlodovsky1925}\cite{kantorovich1931}). For an extensive treatment of the subject, we refer to Ferguson's classical text \cite{ferguson1980}. 
One of the important characteristics of the theory is that uniform approximation is only possible on intervals of length less than $4$, and then only with certain arithmetic conditions on $f$ the number of which increases without bound as the length of the interval increases towards $4$. Kakeya showed (also in \cite{kakeya1914}) that approximating $f \in \mathscr{C}([-2,2])$ arbitrarily well by elements of $\Z[X]$ is impossible unless $f \in \Z[X]$ because metrically, $\Z[X]$ is a discrete subset of $\mathscr{C}([-2,2])$. Indeed,
any non-zero element of $\Z[X]$ has its supremum norm on $[-2,2]$ equal to at least $1$ (and at least $2$ if it is not a constant polynomial). In the other extreme, there are no arithmetic conditions to be met when approximation is sought on closed intervals containing no integers, as shown by Chlodovsky \cite{chlodovsky1925}.

An important case between these two extremes is the interval $[0,1]$ where the characterizing condition reduces to integer valuations at the boundary points $0$ and $1$, as for example, demonstrated by Chlodovsky \cite{chlodovsky1925} and Kantorovich \cite{kantorovich1931}. The corresponding class is
$$\mathscr{C}_\Z([0,1]):=\Big \{f \in \mathscr{C}([0,1]) : f(0),f(1) \in \Z\Big \}.$$
Note that $\mathscr{C}_\Z([0,1])$ is the disjoint union of countably many 	translates of 
$$\mathscr{C}_0([0,1]):=\Big \{f \in \mathscr{C}([0,1]) : f(0)=f(1)=0\Big \}.$$
Indeed, we have
$$\mathscr{C}_\Z([0,1]) = \bigcup_{(m,n) \in \Z^2} \mathscr{C}_0([0,1]) + \ell_{m,n},
$$
where $\ell_{m,n}$ stands for the linear polynomial satisfying $\ell_{m,n}(0)=m$, $\ell_{m,n}(1)=n$.

Let us define $\mathscr{P}_n$ to be the space of polynomials (more precisely, polynomial functions) of degree at most $n$ where, from now on, the domain is understood to be $[0,1]$. If we pick the standard power basis $(\pi_k)_0^n$ for $\mathscr{P}_n$, where $\pi_k(x):=x^k$, then the set of polynomials in $\mathscr{P}_n$ with coefficients in $\Z$ corresponds to the lattice generated by $(\pi_k)_0^n$. We will denote this lattice by $\Pi_n$. Chlodovsky's density result can now be stated as 
\begin{equation}\label{integer_dense}
\bigcup_n \Pi_n, \mbox{ i.e. } \Z[X] \mbox{ is dense in } \mathscr{C}_\Z([0,1]).
\end{equation}

It is natural to ask the rate at which  $d_\infty(f,\Pi_n) \to 0$ given any $f \in \mathscr{C}_\Z([0,1])$. Kantorovich in \cite{kantorovich1931} gave an elementary algorithmic proof of Chlodovsky's theorem that came with a quantitative approximation error bound (see also \cite{ferguson1980} and \cite{lorentz2}). In our notation, Kantorovich's elementary bound states that for all $f \in \mathscr{C}_\Z([0,1])$, 
$$d_\infty(f,\Pi_n):= \inf_{P \in \Pi_n} \|f - P\|_{[0,1]} \leq \|f - B_n\|_{[0,1]} + \frac{1}{n},
$$
where $B_n:=B_n(f)$ is the Bernstein polynomial of $f$ defined by
\begin{equation}\label{def_Bern} 
B_n(x) := \sum_{k=0}^n f\Big(\frac{k}{n}\Big) \binom{n}{k} x^k (1-x)^{n-k}.
\end{equation}
It is well known (see, e.g. \cite{lorentz1}) that 
\begin{equation}\label{Bern_bound}
	\|f - B_n \|_{[0,1]} \leq \frac{5}{4} \omega_f(n^{-1/2}),
\end{equation}
where $\omega_f$ stands for the modulus of continuity of $f$ on $[0,1]$ given by 
$$\omega_f(\delta):=\sup\{|f(x)-f(y)| : x,y \in [0,1], ~|x - y | \leq \delta\},$$
therefore it follows that
$d_\infty(f,\Pi_n) \leq \frac{5}{4} \omega_f(n^{-1/2}) + n^{-1}$.
Kantorovich actually proved a stronger result in \cite{kantorovich1931}, showing that 
$$d_\infty(f,\Pi_n) \leq 2 d_\infty(f,\mathscr{P}_n) + \frac{1}{n}.$$
In other words, up to the additive term $1/n$ and the constant $2$, polynomials with integer coefficients are nearly as expressive as arbitrary polynomials (of degree at most $n$). This was shown by employing the unconstrained polynomial of best approximation of $f$ as a surrogate instead of the Bernstein polynomial (see \cite{ferguson1980}\cite{lorentz1}). Since $d_\infty(f,\mathscr{P}_n) \leq 3\, \omega_f(n^{-1})$ and $\omega_f(\delta) \geq \frac{1}{2} \omega_f(1) \delta$ for all non-constant continuous functions (see e.g. \cite{DL}\cite{lorentz1}) it follows that
$$d_\infty(f,\Pi_n) \leq C_f \omega_f(n^{-1}),
$$
where $C_f \leq 6 + 2/\omega_f(1)$ is a constant that may possibly depend on $f$. (See also \cite[Ch.12]{ferguson1980} as well as \cite{trigub1962}.) 

Finally, as is the case for approximation by polynomials with arbitrary coefficients, the rate of approximation by elements of $\Pi_n$ can be improved under various regularity (e.g. differentiability) assumptions (again, see \cite{ferguson1980}).

In this paper, we are interested in a different form of strengthening the original density theorem of Chlodovsky (and its proof by Kantorovich) by seeking constructive approximations from much smaller lattices. The method of Kantorovich is based on rounding 
%the compound coefficient $f(\frac{k}{n}) \binom{n}{k}$  in the representation \eqref{def_Bern} to a neighboring integer value. It therefore produces 
the coefficients of an approximating polynomial in the lattice $\Pi_n^\circ$ generated by the polynomials $x^k(1-x)^{n-k}$, $k=0,\dots,n$. Clearly, we have $\Pi_n^\circ \subset \Pi_n$, but the simple identity  
$$ x^k = x^k(x+(1-x))^{n-k} = \sum_{l=0}^{n-k} \binom{n{-}k}{l} x^{k+l} (1-x)^{n-(k+l)}$$ 
also shows that $\pi_k \in \Pi_n^\circ$ for all $k=0,\dots,n$, implying that in fact $\Pi_n = \Pi_n^\circ$. Meanwhile, the lattice $\mathscr{B}_n$ generated by the polynomials 
\begin{equation}\label{def_Bern_basis}
p_{n,k}(x) := \binom{n}{k} x^k(1-x)^{n-k}, ~~k=0,\dots,n,
\end{equation}
is a much smaller (coarser) sublattice of $\Pi_n$ since it stretches the $k$th basis polynomial $x^k(1-x)^{n-k}$ by the integer factor $\binom{n}{k}$ which can get exponentially large in $n$. In particular, the fundamental cell of $\mathscr{B}_n$ consists of
$$\prod\limits_{k=0}^n \binom{n}{k}= e^{\frac{1}{2}n^2(1+o(1))}$$ 
distinct elements of $\Pi_n$. 

We call $\mathscr{B}_n$ the {\em Bernstein lattice}. In this paper, we strengthen \eqref{integer_dense} by showing that 
\begin{equation}\label{Bernstein_dense}
	\bigcup_n \mathscr{B}_n \mbox{ is dense in } \mathscr{C}_\Z([0,1]).
\end{equation}
In fact, we prove this result quantitatively by means of our main theorem stated below:
\paragraph{\bf Main Theorem.}
{\em 	For any $f \in \mathscr{C}_\Z([0,1])$ and positive integer $n$, we have
	$$d_\infty(f, \mathscr{B}_n) \leq \frac{5}{4}\omega_f(n^{-1/2}) + \rho(f,n)$$
	where 
	\begin{eqnarray*} \rho(f,n) &:=& \min_{t \in [0,n/2]\cap \Z} \left\{\max\left (\omega_f\Big (\frac{t}{n}\Big ), \frac{1}{2(n+1-2t)}\right ) +
	\frac{1}{\sqrt{2 (t+1)}}\right\}\\
	& \leq & \omega_f\Big (\frac{1}{2} n^{-1/3}\Big) + 2 n^{-1/3}.
	\end{eqnarray*}
	In particular, we have
$$d_\infty(f, \mathscr{B}_n) \leq \frac{9}{4} \omega_f(n^{-1/3}) + 2 n^{-1/3}.$$
}

\paragraph{Note.}
Our motivation and methods were partly developed in the papers \cite{gunturk_li_2023} and \cite{onebitNN}, where we were able to establish pointwise approximation on $(0,1)$ for all $f \in \mathscr{C}([0,1])$, however, our method could not take advantage of the smaller class $\mathscr{C}_\Z([0,1])$ to make the error bound go to $0$ uniformly on $[0,1]$. The present paper achieves the missing uniform approximation on $[0,1]$ by means of a more powerful approximation idea designed for $\mathscr{C}_\Z([0,1])$ with which we are able to control the error bound effectively all the way to the boundary points. However, the price paid for this uniformity is a lower rate of approximation.

\section{Proof of the Main Theorem}

Before we start the proof, let us recall (see, e.g. \cite{DL}\cite{lorentz1}) the following basic polynomial identities which correspond to the central moments, up to degree two, of the Bernstein basis $(p_{n,k})_{k=0}^n$:
\begin{eqnarray}
\sum_{k=0}^n p_{n,k}(x) & ~=~ & 1, \label{moment0} \\
\sum_{k=0}^n (k-nx) p_{n,k}(x) & = & 0, \label{moment1} \\
\sum_{k=0}^n (k-nx)^2 p_{n,k}(x) & = & nx(1-x). \label{moment2} 
\end{eqnarray}
We also note that $p_{n,k} \geq 0$ on $[0,1]$. When interpreting the proof, it is helpful to keep in mind that $p_{n,k}$ peaks at $k/n$, and that the identity \eqref{moment0} corresponds to a smooth partition of unity.  

To begin the proof, let $f \in \mathscr{C}_\Z([0,1])$ and $n$ be a positive integer. For $n=1$, we can approximate $f$ by the linear interpolant $\ell:=f(0)p_{1,0} + f(1) p_{1,1} \in \mathscr{B}_1$ which clearly remains within $\omega_f(1)$ of $f$, so that the stated approximation error bound holds trivially.
Therefore, without loss of generality, let 
$n \geq 2$. 
We will approximate $f$ by a polynomial
$$ Q_n:= \sum_{k=0}^n q_k p_{n,k}
\in \mathscr{B}_n$$ 
which will be constructed with the help of 
two auxiliary polynomials $B_n$ and $P_n$ in $\mathscr{P}_n$ approximating $f$. Note that $B_n$ and $P_n$ will in general not be in $\mathscr{B}_n$, but they will facilitate the approximation and the rounding processes.
The first polynomial $B_n$ is the standard Bernstein polynomial of $f$ of order $n$ defined by \eqref{def_Bern}, i.e.
$$ B_n:= \sum_{k=0}^n f\Big( \frac{k}{n} \Big) p_{n,k}.$$
Once $P_n$ and $Q_n$ are defined, we will simply utilize the triangle inequality
$$ \|f - Q_n \|_{[0,1]} \leq \|f - B_n \|_{[0,1]} + \|B_n - P_n \|_{[0,1]} + \|P_n - Q_n \|_{[0,1]},$$
hence the proof will consist of individual steps bounding each error term.

\par \noindent {\bf 1.} 
For $\|f - B_n\|_{[0,1]}$, it will suffice to use the generic bound \eqref{Bern_bound} stated earlier.

\par \noindent {\bf 2.} The second approximating polynomial $P_n$ will be a perturbation of $B_n$ given by
$$P_n:= \sum_{k=0}^n y_k p_{n,k}$$
where
\begin{equation}
y_k := \left\{ 
\begin{array}{ll}
	f(0), ~&~ 0 \leq k < t_n, \\
	f\big(\frac{k}{n} \big) + \varepsilon_n, ~&~ t_n \leq k \leq n-t_n,\\
	f(1), ~&~ n-t_n < k \leq n,
\end{array}
\right.
\label{def_y_k}
\end{equation}
under the assumption that $t_n$ is a nonnegative integer not exceeding $n/2$ and
\begin{equation}
\varepsilon_n := \frac{1}{n+1-2t_n} \left( 
\left [ \sum_{k=t_n}^{n-t_n} f\Big ( \frac{k}{n} \Big ) \right ]_\Z
- \sum_{k=t_n}^{n-t_n} f\Big ( \frac{k}{n} \Big ) \right).
\label{def_eps_n}
\end{equation}
Here $\big[w\big]_\Z$ stands for any integer such that $\left |w - \big[w\big]_\Z \right | \leq 1/2$, though choosing the closest integer is not critical for our argument.

The significance of the definition \eqref{def_y_k} will become clear when we use the $y_k$ to construct the integer coefficients $q_k$ of the final polynomial $Q_n \in \mathscr{B}_n$ by means of a recurrence relation. For now, we only note that $y_k$ is a uniformly small perturbation of $f(k/n)$. Indeed, 
$$
\left |f\Big(\frac{k}{n} \Big)- y_k \right | \leq
\left\{ 
\begin{array}{ll}
	\omega_f(\frac{k}{n}), ~&~ 0 \leq k < t_n, \\
	|\varepsilon_n|, ~&~ t_n \leq k \leq n-t_n,\\
	\omega_f(\frac{n-k}{n}), ~&~ n-t_n < k \leq n,
\end{array}
\right.
$$
and $|(n+1-2t_n)\varepsilon_n| \leq 1/2$, so that we have 
$$
\max_{0 \leq k \leq n} \left |f\Big(\frac{k}{n} \Big)- y_k \right | \leq \max\left (\omega_f\Big (\frac{t_n}{n}\Big ), \frac{1}{2(n+1-2t_n)}\right ).
$$ 
Therefore, using \eqref{moment0}, it follows that
\begin{equation}
\| B_n - P_n \|_{[0,1]}  \leq \max\left (\omega_f\Big (\frac{t_n}{n}\Big ), \frac{1}{2(n+1-2t_n)}\right ).
\end{equation}

\par \noindent {\bf 3.} We now turn to the construction of $Q_n := \displaystyle \sum_{k=0}^n  q_k p_{n,k} \in \mathscr{B}_n$ and the accompanying bound on $\|P_n - Q_n\|_{[0,1]}$. 
We set $u_{-1}:=0$ and define the $q_k \in \Z$, $k=0,\dots,n$, according to the assignment rule
\begin{eqnarray}
	 q_k &:=& \big [ u_{k-1} + y_k \big ]_\Z,  \label{def_q_k} \\
	 u_k &:=& u_{k-1} + y_k - q_k. \label{recurse}  
\end{eqnarray}
It is apparent from this recurrence that $|u_k| \leq 1/2$ for all $k$. Furthermore, noting that $y_k = f(0) \in \Z$ for $k=0,\dots,t_{n}-1$, it follows recursively that $q_k = y_k = f(0)$ and $u_k=0$ for $k=0,\dots,t_n-1$. Summing \eqref{recurse} for $k=t_n,\dots,n-t_n$ and using that $u_{t_n-1}=0$, we obtain
$$ u_{n-t_n} = \sum_{k=t_n}^{n-t_n} y_k - \sum_{k=t_n}^{n-t_n} q_k.$$
The definition of $y_k$ given in \eqref{def_y_k} and the definition of $\varepsilon_n$ given in \eqref{def_eps_n} together yield
$$
\sum_{k=t_n}^{n-t_n} y_k =  (n+1-2t_n) \varepsilon_n  + 
\sum_{k=t_n}^{n-t_n} f\Big ( \frac{k}{n} \Big ) = 
\left [ \sum_{k=t_n}^{n-t_n} f\Big ( \frac{k}{n} \Big ) \right ]_\Z 
$$
which shows that $u_{n-t_n} \in \Z$. Given that all $u_k$ are bounded by $1/2$ in absolute value, it follows that $u_{n-t_n} = 0$. At this point, noting that $y_k = f(1) \in \Z$ for $k=n-t_n+1,\dots,n$, it follows recursively that $q_k = y_k = f(1)$ and $u_k=0$ for $k=n-t_n+1,\dots,n$. Hence, 
\begin{eqnarray}
P_n - Q_n 
& = & \sum_{k=t_n}^{n-t_n} (y_k - q_k) p_{n,k} \nonumber \\
& = & \sum_{k=t_n}^{n-t_n} (u_k - u_{k-1}) p_{n,k} \nonumber \\
& = & \sum_{k=t_n}^{n-t_n-1} u_k (p_{n,k} - p_{n,k+1}), \label{P_nminusQ_n}
\end{eqnarray}
where in the last equality we have used the fact that $u_{t_n-1} = u_{n-t_n} = 0$. Note that if $t_n = n/2$, then the sum in \eqref{P_nminusQ_n} contains no terms (and $P_n = Q_n$), so from now on we will assume $t_n < n/2$, which guarantees that there is at least one term in this sum.

In order to analyze this representation, we next note the identity (see e.g. \cite{DL},\cite{lorentz1})
\begin{equation}
	\label{p_nk_diff}
x(1-x)\Big(p_{n,k}(x) - p_{n,k+1}(x)\Big) = \Big(\frac{k{+}1}{n{+}1}-x\Big)\,p_{n+1,k+1}(x),
\end{equation}
which holds for all $k \in \Z$, with the natural convention that $p_{n,k} := 0$ whenever $k < 0$ or $k > n$ (since this holds for . Plugging this formula in \eqref{P_nminusQ_n} and noting that $|u_k| \leq 1/2$, we obtain
\begin{eqnarray}
x(1-x)|P_n(x)-Q_n(x)| 
& \leq & \frac{1}{2} \sum_{k=t_n}^{n-t_n-1}
\Big|\frac{k{+}1}{n{+}1}-x\Big|\,p_{n+1,k+1}(x) \nonumber \\
& = & \frac{1}{2} S_{n+1,t_n+1}(x), \label{S_nt_bound_1}
\end{eqnarray}
where
\begin{equation}\label{S_nt_def}
	S_{n,t}(x):=\sum_{k=t}^{n-t}
\Big|\frac{k}{n}-x\Big|\,p_{n,k}(x).
\end{equation}
(Here, note that $t_n < n/2$ implies $t_n+1 \leq (n+1)/2$.)

We claim that for all $x \in [0,1]$ and $0 < t \leq n/2$, we have the bound
\begin{equation} \label{S_nt_bound}
	S_{n,t}(x) \leq \frac{x(1-x)}{\sqrt{t/2}}.
\end{equation}
With this bound, \eqref{S_nt_bound_1} implies (after canceling the term $x(1-x)$ for $0<x<1$ and using continuity at $x=0$ and $x=1$) that 
\begin{equation}
	\|P_n - Q_n\|_{[0,1]} \leq \frac{1}{\sqrt{2 (t_n+1)}}.
\end{equation}

It remains to prove our claim \eqref{S_nt_bound}. Applying
Cauchy-Schwarz to \eqref{S_nt_def} we have
\begin{equation}
S_{n,t}(x) \leq 
\left ( \sum_{k=t}^{n-t} \Big(\frac{k}{n}-x\Big)^2\,p_{n,k}(x) \right)^{1/2}
\left ( \sum_{k=t}^{n-t} p_{n,k}(x) \right)^{1/2}.
\label{S_nt_CS}
\end{equation}
For the first term in this product, \eqref{moment2} implies
\begin{equation}
	\sum_{k=t}^{n-t} \Big(\frac{k}{n}-x\Big)^2\,p_{n,k}(x) 
	\leq \frac{x(1{-}x)}{n}.
	\label{S_nt_term1}
\end{equation}
For the second term, we start by noting that 
$$ k(n-k) \geq t(n-t), ~~k=t,\dots, n-t,
$$
which is a consequence of the fact that the function $g(y):=y(1-y)$ reaches its minimum on the interval $[\tau,1-\tau]$ at $y=\tau$ and $y=1-\tau$, $0\leq \tau \leq 1/2$. Hence 
\begin{eqnarray}
\sum_{k=t}^{n{-}t} p_{n,k}(x) & \leq & \frac{1}{t(n-t)} \sum_{k=t}^{n{-}t} k(n-k) p_{n,k}(x) \nonumber \\
& \leq & \frac{1}{t(n{-}t)} \sum_{k=0}^{n} k(n-k) p_{n,k}(x) \nonumber \\
& = & \frac{n(n{-}1)x(1{-}x)}{t(n{-}t)} \nonumber \\
& < & \frac{2(n{-}1)x(1{-}x)}{t},
\label{S_nt_term2}
\end{eqnarray}
where in the second last step we used the identity
$$
\sum_{k=0}^n k(n-k)p_{n,k}(x) = n(n-1)x(1-x)
$$
which can be easily derived from the moment relations \eqref{moment0}--\eqref{moment2} by means of the identity
$$
k(n-k) = -(k-nx)^2 + n(1-2x)(k-nx) + n^2x(1-x).
$$
Injecting the bounds \eqref{S_nt_term1} and \eqref{S_nt_term2} in \eqref{S_nt_CS} yields the claimed bound in \eqref{S_nt_bound}.

\par \noindent {\bf 4.} 
To complete the proof,
it remains to analyze the sum of these three bounds. With the results of the steps {\bf 1--3}, we have
$$
\|f - Q_n \|_{[0,1]} \leq \frac{5}{4}\omega_f(n^{-1/2}) + \max\left (\omega_f\Big (\frac{t_n}{n}\Big ), \frac{1}{2(n+1-2t_n)}\right ) +
\frac{1}{\sqrt{2 (t_n+1)}}.
$$
So far, the value of $t_n$, $0 \leq t_n \leq n/2$, is arbitrary, noting that $Q_n$ depends only on $f$ and $t_n$ according to the scheme presented in steps {\bf 2} and {\bf 3}. Hence it follows that 
$d_\infty(f, \mathscr{B}_n)$ is bounded above by
$$ \frac{5}{4}\omega_f(n^{-1/2}) + \rho(f,n)
$$
where
$$ \rho(f,n):=  \min_{t \in [0,n/2]\cap \Z} \left\{\max\left (\omega_f\Big (\frac{t}{n}\Big ), \frac{1}{2(n+1-2t)}\right ) +
\frac{1}{\sqrt{2 (t+1)}}\right\}.
$$

We can choose the convenient value $t = t_n := \lfloor \frac{1}{2}n^{2/3} \rfloor$. With this choice, it is straightforward to check that for all $n \geq 2$, we have 
$$
\frac{1}{2(n+1-2t_n)} +
\frac{1}{\sqrt{2 (t_n+1)}} \leq \frac{2}{n^{1/3}}.
$$
Hence it follows that for all $n \geq 2$, we have the bound
$$
d_\infty(f, \mathscr{B}_n) \leq \frac{5}{4} \omega_f(n^{-1/2}) + \omega_f(\frac{1}{2} n^{-1/3}) + 2 n^{-1/3}
$$
which, in particular, is bounded above by the simpler expression 
$$\frac{9}{4}\omega_f(n^{-1/3}) + 2n^{-1/3}.$$
This concludes our proof.
\qed

\paragraph{Notes and Comments.}
\begin{enumerate}
\item Evaluating $\rho(f,n)$ via the optimal choice of $t = t_n$ requires more precise information on $\omega_f$. It can be seen that for Lipschitz functions, the scaling relation $t_n \sim n^{2/3}$ is essentially optimal, hence our choice. However, the situation is different for functions of lower H\"older regularity. If $\omega_f(\delta) \leq C \delta^{1/2}$, then the choice $t_n \sim n^{1/2}$ yields the upper bound $O(n^{-1/4})$ instead of the inferior $O(n^{-1/6})$ that follows from the choice $t_n \sim n^{2/3}$. It can be checked that for H\"older continuous functions with exponent $0 < \alpha \leq 1$ (i.e. $\omega_f(\delta) \leq C \delta^\alpha$), then the optimal exponent $\theta_\alpha$ where $t_n \sim n^{\theta_\alpha}$ is given by
$$
\theta_\alpha = \left \{ 
\begin{array}{ll} 
	\frac{1}{2}, & 0 < \alpha  \leq \frac{1}{2},\\
	\frac{2\alpha}{2\alpha+1}, & \frac{1}{2} \leq \alpha \leq 1,
\end{array} \right.
$$
yielding the bound
$$
d_\infty(f, \mathscr{B}_n) = O\left ( n^{-\min(\frac{\alpha}{2},\frac{\alpha}{2\alpha+1})}\right ).
$$

\item It is natural to ask if the rate of approximation can be improved for differentiable functions, as was the case for the pointwise results in \cite{gunturk_li_2023}. This appears to be nontrivial, if at all possible, the investigation of which we leave for future work.
\item Currently we do not have any expectations on the optimality of our upper bound since we do not know any class of functions with a matching lower bound. As a particular case, it would be interesting to derive a lower bound on the supremum of $d_\infty(f, \mathscr{B}_n)$ over all Lipschitz functions $f \in \mathscr{C}_\Z([0,1])$ with Lipschitz norm at most $1$. A lower bound in this spirit (with an additional constraint that limits the size of the coefficients of the approximation in $\mathscr{B}_n$) was established in \cite{gunturk_li_sampta23}, but without the integer boundary constraint. 

\end{enumerate}


\end{document}